\documentclass[10pt,english]{article}
\usepackage[T1]{fontenc}
\usepackage[latin1]{inputenc}
\usepackage{a4wide}
\usepackage{babel}
\usepackage{psfrag}
\usepackage{latexsym}

\makeatletter

\newtheorem{prop}{Proposition}[section]
\newtheorem{lemma}[prop]{Lemma}
\newtheorem{theorem}[prop]{Theorem}
\newtheorem{defi}[prop]{Definition}
\newtheorem{rem}[prop]{\em Remark\/}
\newtheorem{cor}[prop]{Corollary}
\newtheorem{cond}[prop]{Conditions}
\newtheorem{ex}[prop]{\em Example\/}

\newcommand{\bp}{\begin{prop}}
\newcommand{\bl}{\begin{lemma}}
\newcommand{\bt}{\begin{theorem}}
\newcommand{\bd}{\begin{defi}\rm}
\newcommand{\br}{\begin{rem}\rm}
\newcommand{\be}{\begin{equation}}
\newcommand{\bea}{\begin{eqnarray}}
\newcommand{\bpr}{\begin{proof}}
\newcommand{\bc}{\begin{cor}}
\newcommand{\bco}{\begin{cond}}
\newcommand{\bex}{\begin{ex}\rm}

\newcommand{\ep}{\end{prop}}
\newcommand{\el}{\end{lemma}}
\newcommand{\et}{\end{theorem}}
\newcommand{\ed}{\end{defi}}
\newcommand{\er}{\end{rem}}
\newcommand{\ee}{\end{equation}}
\newcommand{\eea}{\end{eqnarray}}
\newcommand{\epr}{\end{proof}}
\newcommand{\ec}{\end{cor}}
\newcommand{\eco}{\end{cond}}
\newcommand{\eex}{\end{ex}}

\newcommand{\nn}{ \nonumber \\ }
\newcommand{\pr}{{\em Proof.\ }}

\newcommand{\ot}{\otimes}
\newcommand{\ci}{\circ}
\def\op{^{\rm op}}
\newcommand{\M}{\mathcal{M}}

\newcommand{\setc}[1]{\setcounter{equation}{#1}}
\newcommand{\stac}[1]{\stackrel{\ot}{_{_{#1}}}}
\newcommand{\BIM}{{\rm\bf BIM}}

\newcommand{\lb}{\label}

\begin{document}

\large
\title{\bf Internal bialgebroids, entwining structures and corings}
 
\author{\sc Gabriella B\"ohm\\
Research Institute for Particle and Nuclear Physics, Budapest,\\
H-1525 Budapest 114, P.O.B. 49, Hungary\\
E-mail: G.Bohm@rmki.kfki.hu}

\date{}
 
\maketitle
\normalsize 
\vskip 2truecm

\begin{abstract}  
The internal bialgebroid -- in a symmetric monoidal
category with coequalizers -- is defined. The axioms are formulated in
terms of internal 
entwining structures and alternatively, in terms of internal
corings. The Galois property of the coring in 
question is related to the $\times_R$-Hopf algebra property. The language of
entwining structures is used to discuss duality. 
\end{abstract}
\bigskip
\bigskip
\bigskip
\bigskip

\section{Introduction}

Entwining structures have been introduced by T. Brzezi\'nski and
S. Majid \cite{BreMaj} in defining coalgebra principal bundles. The
relevance of entwining structures is raised by the observation of
M. Takeuchi that they provide examples of corings -- i.e. comonoids in
bimodule categories. 

The theory of corings \cite{BreWis} has been
applied among others in the theory of rings and of graded algebras and
also  in non-commutative geometry. One of the most important
applications is, however, its connection with bialgebras. The various
notions of Hopf type modules over a bialgebra can be unified in terms
of comodules over appropriate corings \cite{Bre1,Bre}. In particular,
as it was 
pointed out by R. Wisbauer (\cite{Wis}, Proposition 5.2) an algebra
$(A,\mu_A,\eta_A)$ and a coalgebra $(A,\Delta_A,\epsilon_A)$ over a
commutative ring $k$ (where $A$ is faithful as a $k$-module)
combine into a bialgebra if and only if the $k$-module $A\stac{k} A$
-- with left $A$-action $\mu_{A\stac{k} A}(\Delta_A \stac{k} A \stac{k}
A)$ and right $A$-action $A \stac{k} \mu_A$, coproduct $\Delta_A
\stac{k} A$ and counit $\epsilon_A \stac{k} A$ --
is an $A$-coring. Alternatively, if and only if the algebra and the
coalgebra structures are entwined by the map
$\mu_{A\stac{k}A}(\Delta_A\stac{k} A\stac{k}\eta_A)$.

In the papers \cite{CaenDeGr} and \cite{Wis} entwining structures and
corings have
been generalized by allowing the module structures in the definition
to be non-unital. The main motivation for the introduction of these so
called weak entwining structures and weak corings was to establish
similar connections with weak bialgebras \cite{Nill,BNSz} as one has
between entwining structures, corings and bialgebras.

The motivation of the present paper is similar. In order to make
connection with (internal) 
bialgebroids, we consider internal entwining structures
and corings over monoids in monoidal categories with coequalizers.

The role of corings in the description of Doi-Koppinen modules over a
bialgebroid has been studied already in \cite{BreCaenMil}.
Let ${\mathcal A}=(A,R,s,t,\gamma,\pi)$ be a bialgebroid \cite{Lu,Xu}
or, what is equivalent to it, a $\times_R$-bialgebra \cite{Tak}.
In this notation the algebras $A$ and $R$ over the commutative ring
$k$ are the total and base algebras, respectively, the maps $s:R\to A$
and $t:R\op\to A$ are the source and target maps, $\gamma$ the coproduct
and $\pi$ the counit. Recall that $A$ is an $R$-$R$-bimodule via
$$ r\cdot a\cdot r^\prime\colon = s(r) t(r^\prime) a\qquad
r,r^\prime\in R\quad a\in A.$$
In the paper \cite{BreCaenMil} an $A$-coring has been associated to 
${\mathcal A}$ (and more generally to any Doi-Koppinen datum over
$R$.) It is $A\stac{R} A$ as a $k$-module, the $A-A$-bimodule
structure is given by 
$$ a\cdot (b\stac{R} c)\cdot d= a_{(1)}b \stac{R} a_{(2)}cd \qquad
a,b,c,d\in A$$
(where the usual Sweedler's convention is used and $b\stac{R}c$ is the
image of $b\stac{k}c$ under the canonical projection $A\stac{k}A\to
A\stac{R}A$) and the coproduct and 
the counit are $\gamma\stac{R} A$ and $\pi\stac{R} A$,
respectively. This coring, however, does not come from an entwining
structure.

Our first motivation to consider internal entwining structures is the
aesthetical need to associate an entwining structure (in the bimodule
category ${_R\M_R}$)
to a bialgebroid. This is done in the Example \ref{entwex} below. It
is explained in the Section \ref{relate} how is the internal entwining
structure in the Example \ref{entwex} related to the coring of
\cite{BreCaenMil}. 

Our second motivation comes from the study of internal
bialgebroids. The definition of algebraic structures inside
categories is a very old idea \cite{MacL} with many interesting
results. Let us mention the papers \cite{Par77,Maj,Par97}, for
example. 

Internal bialgebroids are defined in symmetric monoidal categories
with coequalizers. They consist of the data $(A,R,s,t,\gamma,\pi)$
where $A$ and $R$ are monoids in $\M$ and $s,t,\gamma$ and $\pi$ are
morphisms in $\M$. The axioms are naively copied from the case when
$\M$ is the category of modules over a commutative ring $k$.

Bialgebroids have been characterized by B. Day and R. Street in purely
categorical terms \cite{DS}. They consider two objects $R$ and $A$ in
a monoidal bicategory, a pseudo-monoid structure on $A$ together with
a certain strong monoidal morphism from this pseudo-monoid to the
canonical  pseudo-monoid associated to $R$. In the particular case of
the monoidal bicategory of bimodules -- i.e. the one with 0-cells the
$k$-algebras, 1-cells the bimodules and 2-cells the bimodule
morphisms, --  this definition recovers the one of the bialgebroid.
As a support of our Definition \ref{bgd} of internal
bialgebroids we prove an equivalence with the description of \cite{DS}
in the case of the monoidal bicategory of internal bimodules.

We give a characterization of internal bialgebroids also in
terms of internal entwining structures and alternatively in terms of
internal corings. An (internal) coring
associated to an (internal) bialgebroid is shown to be Galois (in the
sense of 
Definition \ref{Galois}) if and only if the bialgebroid is a
$\times_R$-Hopf algebra in the sense of \cite{Schau}.

In the case of a bialgebroid 
${\mathcal A}=(A,R,s,t,\gamma,\pi)$ over the commutative ring $k$ one of the
$R$-module  
structures of $A$ was assumed 
to be finitely generated projective in \cite{KadSzlach}. Under this
assumption the $R$-dual 
was shown to have an opposite bialgebroid structure (which was called a
right bialgebroid there).

Our formulation of the axioms in terms of an entwining structure provides a
natural framework for the study of the duality of internal bialgebroids, as
the entwining structure is a self-dual notion. Generalizing the result of
\cite{KadSzlach} we prove that if both objects
appearing in the entwining structure possess right duals, then the dual
entwining structure determines a dual internal (right) bialgebroid.  

\section{Preliminaries}
\lb{preli}
\setc{0}

In this section we fix our notations and conventions.

Throughout the paper $k$ is a commutative ring and $\M_k$ denotes the
symmetric monoidal category of $k$ modules. 

The monoidal categories of the paper are {not} required to be
strict but -- relying on coherence -- we do not denote the
reassociator and unit isomorphisms.

Let $(\M,\Box,U)$ be a monoidal category which
possesses coequalizers and in which the monoidal product
preserves the coequalizers. (In this situation we say shortly that $\M$
is a monoidal category with coequalizers. A typical example
of such categories is $\M_k$.) 
Let us construct the bicategory $\BIM(\M)$. The 0-cells in
$\BIM(\M)$ are the monoids in $\M$, the 1-cells the bimodules in $\M$
and the 2-cells the internal bimodule morphisms. We use bold face
letters for the cells in $\BIM(\M)$ and the same Roman letters for
their $\M$-object part. For example we write ${\bf
A}=(A,\mu_A,\eta_A)$ for a monoid and ${\bf M}=(M,\lambda_M,\rho_M)$ for
a bimodule in $\M$.

\begin{figure}[h]
\psfrag{1}{\Large$M\Box B\Box N$}
\psfrag{2}{\Large$M\Box N$}
\psfrag{3}{\Large$M\stac{B} N$}
\psfrag{4}{$\rho_M\Box N$}
\psfrag{5}{$M\Box\lambda_N$}
\psfrag{6}{$\sqcap(M,N)$}
\begin{center}
{\resizebox*{8cm}{!}{\includegraphics{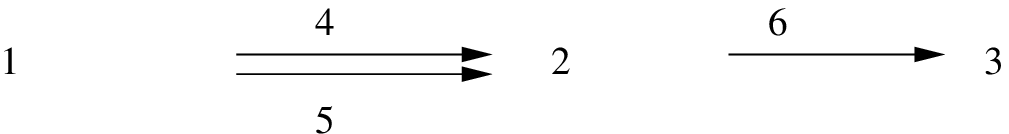}}
}
\end{center}
\caption{The coequalizer diagram used to define the module tensor
product}
\label{modprod}
\end{figure}

The vertical product in $\BIM(\M)$ is the composition $\ci$ in $\M$ and the
horizontal product is the module tensor product -- constructed with the
help of the coequalizers as follows: Let ${\bf A},{\bf B},{\bf C}$ be
monoids in $\M$ and let ${\bf M}=(M,\lambda_M,\rho_M)$ and ${\bf
N}=(N,\lambda_N,\rho_N)$ be ${\bf A}$-${\bf B}$ and ${\bf B}$-${\bf C}$
bimodules, respectively. The object $M\stac{B}N$ is the object part of
the coequalizer of the parallel morphisms $\rho_M\Box N$ and $M\Box
\lambda_N$ as on Figure \ref{modprod}.
We use the notations $\sqcap(\ ,\ )$ and  $\stac{B}$ as on Figure
\ref{modprod} throughout the
paper. The ${\bf A}$-${\bf C}$-bimodule ${\bf M\stac{B} N}=(M\stac{B}
N,\lambda_{M\stac{B} N},\rho_{M\stac{B} N})$ consists of the unique
morphisms $\lambda_{M\stac{B} N}$ and $\rho_{M\stac{B} N}$ in $\M$
which satisfy
\bea \lambda_{M\stac{B} N}
\ci \left(A\Box \sqcap(M,N)\right)&=&
\sqcap(M,N)\ci \left(\lambda_M\Box N\right)\quad {\textrm {and}}\\
\rho_{M\stac{B} N}\ci\left(\sqcap(M,N)\Box C\right)&=&
\sqcap(M,N)\ci \left(M\Box\rho_N\right), 
\nonumber\eea
respectively. The ${\bf B}$-module tensor product of bimodule
morphisms ${\bf \xi}:{\bf M}\to {\bf M^{\prime}}$ and ${\bf \zeta}:
{\bf N}\to {\bf N^{\prime}}$ is the unique morphism $\xi\stac{B}\zeta$
in ${\M}$ which satisfies
$$ (\xi\stac{B}\zeta)\ci\sqcap(M,N)=\sqcap(M^{\prime},N^{\prime}) 
\ci (\xi\Box\zeta). $$
The morphisms $\lambda_{M\stac{B} N}$, $\rho_{M\stac{B} N}$ and
$\xi\stac{B}\zeta$ have been constructed using the universality of the
coequalizer -- what will be done often in the paper. Recall that such
definitions ``$g$ is the unique morphism for which $g\ci \sqcap(M,N)
=f$ '' make sense provided $f\ci (\rho_M\Box
N)=f\ci(M\Box\lambda_N)$. 

Recall also that the coequalizer $\sqcap(M,N)$ is an epimorphism
hence we can transform statements on 1-cells in $\BIM(\M)$ to
equivalent statements on morphisms in $\M$ by composing with it on the
right. 

About the construction of the coherent reassociator and unit
isomorphisms in $\BIM(\M)$ consult the Appendix in \cite{FrobD2}.

Relying on coherence again, we do not denote the reassociator
isomorphism in $\BIM(\M)$ and identify the isomorphic objects
$A\stac{A} M$ and $M\stac{B}B$ with $M$.

If the category $\M$ is also symmetric then the bicategory $\BIM(\M)$
is also monoidal.
The monoidal product is the monoidal product $\Box$ of
$\M$ on the objects with the obvious inducement on the morphisms
\bea {\bf A}\Box {\bf B}&=& (A\Box B,
(\mu_A\Box \mu_B)\ci(A\Box \Sigma_{B,A}\Box B),
\eta_A\Box\eta_B)\nn
{\bf M}\Box {\bf N} &=& (M\Box N,
(\lambda_M\Box \lambda_N)\ci(A\Box \Sigma_{B,M}\Box N),
(\rho_M\Box \rho_N)\ci (M\Box \Sigma_{N,C}\Box D) )
\nonumber \eea
for monoids ${\bf A},{\bf B},{\bf C},{\bf D}$, the ${\bf A}$-${\bf
C}$-bimodule ${\bf M}$ and ${\bf B}$-${\bf D}$-bimodule ${\bf N}$.

\section{Entwining structures in monoidal categories}
\setc{0}

The classical notion of entwining structure \cite{BreMaj} consists of
an algebra $S$ and 
a coalgebra $L$ over a commutative ring $k$ and a linear map
$\psi:S\stac{k} L\to L\stac{k} S$ 
relating the algebra and the coalgebra structures. In what follows we
replace the algebra $S$ and the coalgebra $L$ by a monoid ${\bf S}$
and a comonoid ${\bf L}$ in a monoidal category $(\M,\Box,U)$.
\bd \lb{entw}
A {\em left entwining structure} in a monoidal category $(\M,\Box,U)$
consists of a monoid ${\bf S}=(S,\mu,\eta)$ a comonoid ${\bf
L}=(L,\gamma,\pi)$ and a morphism $\psi:S\Box L \to L\Box S$ in $\M$
satisfying the conditions
\bea 
\lb{entwi} &&\psi\ci(\eta \Box L)=L\Box \eta\\
\lb{entwii} &&(\pi\Box S)\ci \psi= S\Box \pi\\
\lb{entwiii} &&(L\Box \mu)\ci (\psi\Box S)\ci (S\Box \psi) = \psi\ci
(\mu\Box L) \\
\lb{entwiv} &&(\gamma \Box S)\ci \psi=(L\Box \psi)\ci (\psi\Box L)\ci
(S\Box \gamma).
\eea
\ed
\bex {\em  Mixed distributive laws.}
A left entwining structure in the monoidal category of endofunctors on a
category ${\cal C}$ is the same as a mixed distributive law \cite{Beck}
of a monad over a comonad on ${\cal C}$.
\eex 
\bex {\em The internal entwining structure associated to a
bialgebroid.} 
\lb{entwex}
Let $(A,R,s,t,\gamma,\pi)$ be a bialgebroid. The bimodule 
$$ S:\qquad r\cdot a\cdot r^{\prime}=s(r)as(r^{\prime}) \qquad
r,r^{\prime}\in R,\ a\in A $$
with the multiplication 
$$ \mu:S\stac{R} S\to S\qquad a\stac{R} b \mapsto ab $$
and unit $s:R\to A$ form a monoid in the monoidal category
$({_R\M_R},\stac{R},R)$ of $R$-$R$-bimodules.

The bimodule 
$$ L:\qquad r\cdot a\cdot r^{\prime}=s(r)t(r^{\prime})a\qquad
r,r^{\prime}\in R,\ a\in A $$
with the comultiplication $\gamma$ and counit $\pi$ form a comonoid in
$({_R\M_R},\stac{R},R)$.

Introducing the bimodule morphism
$$\psi:S \stac{R} L\to  L\stac{R} S\qquad a\stac{R} b \mapsto 
a_{(1)}b\stac{R} a_{(2)} $$
(where the usual Sweedler's convention has been used) we have a left
entwining structure in ${_R\M_R}$.
\eex
\br \lb{entwdual}
The notion of the internal entwining structure is self-dual in the following
sense. Let the monoid ${\bf S}=(S,\mu,\eta)$, the comonoid ${\bf
L}=(L,\gamma,\pi)$ and the morphism $\psi:S\Box L\to L\Box S$ form a left
entwining structure in the monoidal category $(\M,\Box,U)$. Suppose
furthermore that both objects $S$ and $L$ possess right duals $S^r$ and $L^r$
in $\M$. Then the monoid ${\bf L}^r=(L^r,\gamma^r,\pi^r)$, the comonoid ${\bf
S}^r=(S^r,\mu^r,\eta^r)$ and the morphism $\psi^r:S^r\Box L^r\to L^r\Box
S^r$ form a {\em right entwining structure} that is they satisfy
\bea
\lb{rentwi}&&\psi^r\ci (S^r\Box \pi^r)=\pi^r\Box S^r\\
\lb{rentwii}&&(L^r\Box \eta^r)\ci \psi^r=\eta^r\Box L^r\\
\lb{rentwiii}&&(\psi^r\Box S^r)\ci (S^r\Box \psi^r)\ci (\mu^r\Box L^r)=
(L^r\Box \mu^r)\ci \psi^r\\
\lb{rentwiv}&&\psi^r\ci (S^r\Box \gamma^r)=(\gamma^r\Box S^r)\ci 
(L^r\Box \psi^r)\ci(\psi^r\Box L^r).
\eea
\er
\section{Corings in monoidal categories with coequalizers}
\setc{0}

The classical corings \cite{Swe} are comonoids in
$({_A\M_A},\stac{A},A)$, the monoidal
category of bimodules over a $k$-algebra $A$. In what follows we are
going to replace $A$ with a monoid ${\bf A}$ in a monoidal
category $(\M,\Box,U)$. In order to define module tensor
product (over ${\bf A}$) we need the further assumption that
$(\M,\Box,U)$ is a monoidal category with
coequalizers. As it is explained in the Section \ref{preli}, under this
assumption the category ${_{\bf A} \M _{\bf A}}$ of ${\bf A}$-${\bf
A}$-bimodules in $\M$ is a (lax) monoidal category with monoidal unit
${\bf A}$ and monoidal product $\stac{A}$.
\bd Let $(\M,\Box,U)$ be a monoidal category with
coequalizers and let ${\bf A}$ be a monoid in $\M$. {\em An ${\bf
A}$-coring in ${\M}$} is a comonoid in ${_{\bf A} \M _{\bf A}}$. The
morphisms of ${\bf A}$-corings in $\M$ are the comonoidal morphisms in
${_{\bf A} \M _{\bf A}}$.
\ed
We use the notation ${\bf C}=(C,\Lambda,{\cal P},\Delta,\epsilon)$ for
an ${\bf A}$-coring where $(C,\Lambda,{\cal P})$ is an ${\bf A}$-${\bf
A}$-bimodule in $\M$, the $\Delta$ is the coproduct and $\epsilon$ is the
counit. 
\bex A coring over a $k$-algebra $A$ is a coring in the 
monoidal category of $k$-modules.
\eex
\bex Any comonoid in a monoidal category with coequalizers
is a coring over the monoidal unit (with its obvious monoid
structure).
\eex
\bex \lb{Swe}
{\em The internal Sweedler coring.} Let $(\M,\Box,U)$ be
a monoidal category with coequalizers and ${\bf
A}=(A,\mu,\eta)$ a monoid in $\M$. The ${\bf A}$-${\bf
A}$-bimodule $(A\Box A, 
\mu\Box A, A\Box \mu)$ equipped with the coproduct $A\Box \eta \Box A:
A\Box A \to (A\Box A)\stac{A} (A\Box A)\simeq A\Box A\Box A$ and
counit $\mu:A\Box A\to A$ is an ${\bf A}$ coring in $\M$.

In particular let $\iota:R\to A$ be an extension of $k$-algebras. Then
the canonical $R$-$R$-bimodule 
$$ r\cdot a \cdot r^{\prime}=\iota(r)a\iota(r^{\prime})\qquad
r,r^{\prime}\in R\quad  a\in A$$
with the multiplication $\mu:A\stac{R}A\to A$, $a\stac{R}
a^{\prime}\mapsto a a ^{\prime}$ and unit $\iota$ is a monoid ${\bf
A}$ in
${_R\M_R}$ and $A\stac{R} A$ has an ${\bf A}$-coring structure.
\eex
\bex \lb{entw-coring}
{\em Corings from entwining structures.} Generalizing the
observation of Takeuchi (\cite{BreWis}, Proposition 32.6) one can
construct examples of corings in monoidal categories $(\M,\Box,U)$ 
with coequalizers from entwining structures in it. As
a matter of fact let ${\bf S}=(S,\mu,\eta)$ be a monoid, ${\bf
L}=(L,\gamma,\pi)$ a comonoid and $\psi:S\Box L\to L\Box S$ a
morphism in $\M$. Introduce the morphism $\Lambda\colon = (L\Box
\mu)\ci(\psi\Box S): S\Box L\Box S \to L\Box S$ in $\M$. Then the
datum
$$ \left( L\Box S,\Lambda,L\Box \mu,\gamma\Box S,\pi\Box S\right) $$
is an ${\bf S}$-coring in $\M$ if and only if the triple $({\bf
S},{\bf L},\psi)$ is a left entwining structure in $\M$.
\eex
A most classical example of corings (in ${\M_k}$) is the Sweedler coring
with object part $A\stac{R}A$ for an extension of $k$-algebras $R\to
A$. Our Example \ref{Swe} leads, however, to an internal coring in
${_R\M_R}$. 

In \cite{BreCaenMil} a coring (in $\M_k$) has been associated to a
bialgebroid. Applying however the construction of Example
\ref{entw-coring} to the internal entwining structure in the Example
\ref{entwex} we obtain an internal coring in ${_R\M_R}$.

\begin{figure}[h]
\psfrag{1}{\Large$F(M)\Diamond F(A)\Diamond F(N)$}
\psfrag{2}{\Large$F(M)\Diamond F(N)$}
\psfrag{3}{\Large$F(M\stac{A} N)$}
\psfrag{4}{$F(\rho_M)\ci F_2\Diamond F(N)$}
\psfrag{5}{$F(M)\Diamond F(\lambda_N)\ci F_2$}
\psfrag{6}{$F(\sqcap(M,N))\ci F_2$}
\begin{flushleft}
{\resizebox*{13.5cm}{!}{\includegraphics{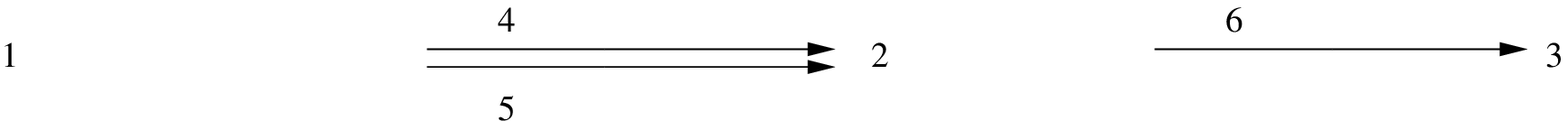}}
}
\end{flushleft}
\caption{The coequalizer diagram in ${\cal N}$}
\label{Fpres}
\end{figure}

The following proposition shows that these corings in $\M_k$ can
be obtained from the ones in ${_R\M_R}$ by applying the
forgetful functor ${_R\M_R}\to \M_k$. 
\bp \lb{forget}
Let $(\M,\Box,U)$ and $({\cal N},\Diamond,V)$ be monoidal categories with
coequalizers. Let $(F,F_0,F_2): {\bf{\cal M}}\to {\bf{\cal N}} $ be 
a (lax) monoidal functor which preserves the module tensor products i.e. 
such that the diagram on Figure \ref{Fpres} is a coequalizer 
diagram in ${\cal N}$ for any monoid ${\bf A}$ and bimodules ${\bf M}$
and ${\bf N}$ in $\M$. Then the functor $F$ preserves the corings.

If furthermore $F$ is faithful and $F_2$ is epi then $F$ also reflects
the corings.
\ep
\br The module tensor product preserving property of a monoidal functor 
$F:{\cal M}\to {\cal N}$ is equivalent to the requirement that the composite
functor $F G^A:{_{\bf A}\M_{\bf A}}\to {\cal N}$ is {\em essentially strong
monoidal} \cite{Szlach3} for any monoid ${\bf A}$ in $\M$, where $G^A$ is the
forgetful functor ${_{\bf A}\M_{\bf A}}\to \M$.
\er
{\em Proof of Proposition \ref{forget}.} 
Any monoidal functor $(F,F_0,F_2)$ preserves monoids, bimodules 
and bimodule morphisms. As a matter of fact for a monoid ${\bf A}=(A,
\mu,\eta)$, an ${\bf A}$-${\bf A}$ bimodule ${\bf M}=(M,\lambda,\rho)$
and  ${\bf A}$-${\bf A}$ bimodule morphism $\xi$ in $\M$ the $F({\bf
A})=(F(A), F(\mu)\ci F_2,F(\eta)\ci F_0)$ is a monoid, $F({\bf
M})=(F(M), F(\lambda)\ci F_2, F(\rho)\ci F_2)$ is an  $F({\bf
A})$-$F({\bf A})$ -bimodule and $F(\xi)$ is an $F({\bf A})$-$F({\bf
A})$-bimodule morphism in ${\cal N}$.

Let us define the functor ${\tilde F}:{_{\bf A}\M_{\bf A}}\to 
{_{F({\bf A})} {\cal N}_{F({\bf A})}}$ as
$$ {\bf M}\mapsto F({\bf M})\qquad \xi\mapsto F(\xi) $$
on the objects and on the morphisms, respectively. Then ${\tilde F}$
is a (lax) monoidal functor with ${\tilde F}_0=F(A)$ and ${\tilde
F}_2$ the unique morphism in ${\cal N}$ for which
$$ {\tilde F}_2 \ci \sqcap(F(M),F(N))=F(\sqcap(M,N))\ci F_2. $$

Now by the assumption that $F$ preserves the module tensor product
${\tilde F}_2=F(M\stac{A} N)$ hence ${\tilde F}$ is strict 
monoidal. Since strong monoidal functors preserve comonoids this
proves that for an ${\bf A}$-coring $(C,\Lambda,{\cal P},\Delta,\epsilon)$
in $\M$ the $(F(C),F(\Lambda)\ci F_2,F({\cal P})\ci
F_2,F(\Delta),F(\epsilon))$ is an $F({\bf A})$-coring in ${\cal N}$. 

If $F$ is also faithful and $F_2$ is epi then $F$ reflects the
monoids, bimodules and bimodule morphisms. Since ${\tilde F}$ is
strong monoidal, it reflects the comonoids. Hence if
$(F(C),F(\Lambda)\ci F_2,F({\cal P})\ci F_2,F(\Delta),F(\epsilon))$ is an
$F({\bf A})$-coring in ${\cal N}$ for some monoid ${\bf A}$, object
$C$ and morphisms $\Lambda,{\cal P},\Delta$ and $\epsilon$ in $\M$ then
$(C,\Lambda,{\cal P},\Delta,\epsilon)$ is an ${\bf A}$-coring in $\M$.
\hfill $\Box$

\bigskip
The following definition of the Galois property of internal corings is
a straightforward generalization of the definition in \cite{Bre}. In
what follows $(\M,\Box,U)$ is a monoidal category
with coequalizers, ${\bf A}=(A,\mu,\eta)$ is a monoid and
${\bf C}=(C,\Lambda,{\cal P},\Delta,\epsilon)$ is an ${\bf A}$-coring
in $\M$. 
\bd A {\em right ${\bf C}$-comodule in $\M$} is a pair $({\bf
M},\tau_M)$ where ${\bf M}$ is a right ${\bf A}$-module  and
$\tau_M:M\to M\stac{A} C$ is a right ${\bf A}$-module morphism in $\M$
satisfying 
\bea \lb{comcp} &&(\tau_M\stac{A} C)\ci \tau_M=(M\stac{A} \Delta)\ci
\tau_M \\
\lb{comcu}&& M=(M\stac{A}\epsilon)\ci \tau_M.
\eea
A right ${\bf C}$ comodule morphism $({\bf M},\tau_M)\to ({\bf
N},\tau_N)$ is a right ${\bf A}$-module morphism $\phi:{\bf M}\to {\bf
N}$ satisfying 
$$ \tau_N\ci \phi=(\phi\stac{A} C)\ci \tau_M. $$
%The category of right ${\bf C}$-comodules in $\M$ is denoted by
%$\M^C$.
\ed
\bp The right ${\bf A}$-module $(A,\mu)$ in $\M$ can be equipped with
a right ${\bf C}$-comodule structure if and only if there exists a
morphism $g:U\to C$ satisfying the conditions
\bea \lb{gcp} && \Delta\ci g=\sqcap(C,C)\ci (g\Box g)\\
\lb{gcu}&& \epsilon \ci g=\eta.
\eea
The morphism $g$ is called a {\em group-like morphism} for ${\bf C}$.
\ep
{\em Proof.}
Suppose that $\tau_A:A\to C$ is a right ${\bf C}$ coaction and set
$g=\tau_A\ci \eta_A$. It satisfies (\ref{gcp}) by (\ref{comcp}) and
(\ref{gcu}) by (\ref{comcu}). Conversely, suppose that there exists a
group-like morphism $g$ for ${\bf C}$ and set $\tau_A={\cal P}\ci
(g\Box A):A\to C$. It is obviously a right ${\bf A}$-module
morphism. The identity (\ref{comcp}) follows from (\ref{gcp}) and
(\ref{comcu}) from (\ref{gcu}).
\hfill $\Box$

\begin{figure}[h]
\psfrag{3}{\Large$M^{\bf C}_g$}
\psfrag{1}{\Large$M$}
\psfrag{2}{\Large$M\stac{A} C$}
\psfrag{6}{$\iota_M$}
\psfrag{4}{$\tau_M$}
\psfrag{5}{$\sqcap(M,C)\ci(M\Box g)$}
\begin{center}
{\resizebox*{6.5cm}{!}{\includegraphics{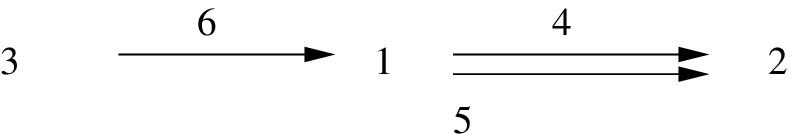}}
}
\end{center}
\caption{The definition of the coinvariant subobject}
\label{inv}
\end{figure}

\bd \lb{invdef}
Let ${\bf C}$ be an ${\bf A}$-coring with group-like
morphism $g$ and let $(M,\tau_M)$ be a right ${\bf C}$-comodule in
$\M$. 
If there exists the equalizer of the parallel morphisms $\tau_M$ and
$\sqcap(M,C)\ci(M\Box g)$ in $\M$ as on the Figure \ref{inv}, then the
object part $M^{\bf C}_g$ of the equalizer is the {\em coinvariant
subobject} of $M$ with respect to $g$.
In particular, if there exists the equalizer of the parallel morphisms
${\cal P}\ci (g\Box A)$ and $\Lambda \ci (A\Box g)$ in $\M$ then its
object part is the coinvariant subobject of $A$ with respect to $g$.
\ed
Let ${\bf C}$ be an ${\bf A}$-coring in $\M$ with group-like
morphism $g$ and corresponding coinvariant subobject $B$ of $A$. The
object $B$ -- if exists -- has  a monoid structure with 
multiplication $\mu_B$ and unit $\eta_B$ which are the unique
morphisms in $\M$ which satisfy 
\bea &&\iota_A\ci \mu_B=\mu\ci (\iota_A\Box \iota_A)\quad {\textrm {and}}\nn
&&\iota_A\ci \eta_B=\eta.
\nonumber\eea
Denote the ${\bf A}$-${\bf B}$-bimodule $(A,\mu,\mu\ci(A\Box \iota_A)
)$ by $\iota_*$ and the ${\bf B}$-${\bf A}$-bimodule $(A,\mu\ci(\iota_A
\Box A) ,\mu)$ by $\iota^*$. The ${\bf A}$-${\bf A}$-bimodule
$\iota_*\stac{B} \iota^*$ is completed to an ${\bf A}$-coring ${\bf
C}_B$ with coproduct $\Delta_B=\iota_*\stac{B} \iota_A\stac{B}
\iota^*$ and counit $\epsilon_B$ which is the
unique morphisms in $\M$ for which
$$\epsilon_B\ci \sqcap(\iota_*,\iota^*)=\mu.$$
\bd \lb{Galois}
The ${\bf A}$-coring ${\bf C}$ with group-like morphism $g$ and
corresponding coinvariant subobject $B$ of $A$ is {\em Galois} (with
respect to $B$) if there
exists an isomorphism of ${\bf A}$-corings $\kappa:{\bf C}_B\to {\bf
C}$ such that 
\be \lb{kappaprop}\kappa\ci\sqcap(\iota_*,\iota^*)\ci (\eta\Box
\eta)=g.\ee 
\ed

\bigskip
%\noindent
The only ${\bf A}$-${\bf A}$-bimodule morphism
$\kappa:\iota_*\stac{B}\iota^*\to (C,\Lambda,{\cal P})$ which is
subject to (\ref{kappaprop}) satisfies
\be\lb{kappaform}
\kappa\ci\sqcap(\iota_*,\iota^*)=\Lambda\ci[A\Box {\cal P}\ci(g\Box
A)].\ee
Since a morphism satisfying (\ref{kappaform}) is a morphism of ${\bf
A}$-corings ${\bf C_B}\to {\bf C}$, 
the Galois property of ${\bf C}$ w.r.t. $B$ is equivalent to the
isomorphism property of the morphism $\kappa$ defined by
(\ref{kappaform}). 

\section{Bialgebroids in symmetric monoidal categories with
coequalizers}
\setc{0}

Throughout the section $(\M,\Box,U,\Sigma)$ is a symmetric monoidal
category with coequalizers.
\bd \lb{bgd}
A {\em left bialgebroid in $\M$} consists of the data 
\begin{itemize}
\item monoids ${\bf R}=(R,\mu_R,\eta_R)$ and ${\bf
A}=(A,\mu_A,\eta_A)$ in $\M$
\item monoidal morphisms $s:{\bf R}\to {\bf A}$ and $t:{\bf R}^{op}\to
{\bf A}$ satisfying
$$ \mu_A\ci(s\Box t)=\mu_A\ci \Sigma_{A,A}\ci(s\Box t)$$
\item a comonoid ${\bf L}=(L,\gamma,\pi)$ in ${_{\bf R}\M_{\bf R}}$
where ${ L}$ is the internal ${\bf R}$-${\bf R}$-bimodule
\be\lb{L} (A,\mu_A\ci (s\Box A), \mu_A\ci \Sigma_{A,A}\ci (A\Box t)
)\ee
\end{itemize}
subject to the axioms
\bea 
\lb{bgdi} &&\rho_{M_A}\ci (\gamma\Box t\Box \eta_A)=\rho_{M_A}\ci
(\gamma\Box\eta_A\Box s)\\
\lb{bgdii}&&\gamma\ci \eta_A=\sqcap(L,L)\ci(\eta_A\Box \eta_A)\\
\lb{bgdiii}&&\gamma\ci \mu_A=\lambda_{M_A}\ci(A\Box \gamma)\\
\lb{bgdiv}&&\pi\ci \eta_A=\eta_R\\
\lb{bgdv}&&\pi\ci\mu_A\ci (A\Box s\ci \pi)=\pi\ci\mu_A=\pi\ci\mu_A\ci
(A\Box t\ci \pi)
\eea 
where we introduced the ${\bf A}$-${\bf A}\Box {\bf A}$-bimodule ${\bf
{M_A}}=( L\stac{R}L,\lambda_{M_A},\rho_{M_A})$ with the unique
morphisms $\lambda_{M_A}$ and $\rho_{M_A}$ which satisfy
\bea \lb{rho}&&\rho_{M_A}\ci (\sqcap(L,L)\Box A\Box A)=\sqcap(L,L)\ci
\mu_{A\Box A}\\
\lb{lambda}&&\lambda_{M_A}\ci(A\Box \sqcap(L,L))=\rho_{M_A}\ci
(\gamma\Box A\Box A).
\eea
(The definition (\ref{lambda}) makes sense by (\ref{bgdi}).) We use
the notation ${\cal A}=({\bf A},{\bf R},s,t,\gamma,\pi)$ for a
bialgebroid.
\ed
Applying the Definition \ref{bgd} in the category $\M=\M_k$ of modules
over a commutative ring $k$ we recover the usual definition
\cite{Tak,Lu,Xu,Szlach2} of the bialgebroid.

The morphisms of left bialgebroids over a fixed monoid ${\bf R}$ in $\M$
are introduced as follows. The bialgebroid morphisms from ${\cal A}=({\bf
A},{\bf R},s,t,\gamma,\pi)$ to ${\cal A}^{\prime}=({\bf A}^{\prime},
{\bf R},s^{\prime},t^{\prime},\gamma^{\prime},\pi^{\prime})$ 
are the morphisms $\Phi:A\to A^{\prime}$ in $\M$ which satisfy
\bea \lb{bgdmori}&& \Phi\ {\textrm {is\ a\ monoidal\ morphism\ }}{\bf
A}\to {\bf A}^{\prime}\ {\textrm{in}}\ \M\\
\lb{bgdmorii}&&\Phi\ {\textrm {is\ a\ comonoidal\ morphism\ }}{\bf
L}\to {\bf L}^{\prime}\ {\textrm{in}}\ {_{\bf R}\M_{\bf R}}\ .
\eea 
The category with objects the left bialgebroids in $\M$ over the base
${\bf R}$ and morphisms defined this way is denoted by ${\bf
Bgd_{R}}(\M)$. 

\bigskip

\noindent
In the paper \cite{DS} it has been proven that a bialgebroid in
$\M_k$ is equivalent to the data 
\begin{itemize}
\item two $k$-algebras ${\bf R}$ and ${\bf A}$ 
\item a $k$- algebra map $\chi:{\bf R^e}\to {\bf A}$ 
\item a pseudo-monoid $({\bf A},{\bf M},{\bf J})$ in $\BIM(\M_k)$
\end{itemize}
such that the 
${\bf R^e}$-${\bf
A}$-bimodule induced by $\chi$ -- that is ${\bf
  \chi^*}=(A,\mu_A\ci(\chi\Box A),\mu_A)$ -- is a strong 
monoidal morphism from $({\bf A},{\bf M},{\bf J})$ to the canonical
pseudo-monoid $({\bf R^e},{\bf m},{\bf j})$. (The symbol ${\bf R^e}$
stands for the enveloping algebra ${\bf R\Box R}^{op}$.)
In the following we prove
an analogous result on bialgebroids in arbitrary symmetric monoidal
category $\M$ with coequalizers.

Recall that a pseudo-monoid in $\BIM(\M)$ consists of a monoid ${\bf
A}$, an ${\bf A}$-${\bf A}\Box {\bf A}$-bimodule ${\bf M}$ and a left
${\bf A}$-module ${\bf J}$ together with bimodule isomorphisms
\bea 
&&l_A:{\bf M}\stac{{\bf A}\Box {\bf A}} ({\bf J}\Box {\bf A})\to {\bf A}\nn
&&r_A:{\bf M}\stac{{\bf A}\Box {\bf A}} ({\bf A}\Box {\bf J})\to {\bf
A}\nn
&&a_A:{\bf M}\stac{{\bf A}\Box {\bf A}} ({\bf M}\Box {\bf A})\to 
{\bf M}\stac{{\bf A}\Box {\bf A}} ({\bf A}\Box {\bf M})
\nonumber\eea 
satisfying the pentagon and triangle identities.

Since ${\bf R^{op}}=(R,\mu_R\ci \Sigma_{R,R},\eta_R)$ is the two-sided
pseudo-dual of ${\bf R}$ in $\BIM(\M)$, the monoid ${\bf R^e}\equiv 
{\bf R}\Box {\bf R^{op}}$ is the monoid part of the canonical pseudo-monoid
with bimodules
\bea \lb{m} &&
{\bf m}=\left(R\Box R\Box R, 
\lambda_m=[\mu_R\Box (R\Box \mu_R)]\ci (R\Box\Sigma_{R, R^{\Box 3}}),\right.\nn
&&\left.\ \qquad
\rho_m=[\mu_R\Box \mu_R\ci(\mu_R\Box R)\Box \mu_R]\ci 
(R^{\Box 4}\Box \Sigma_{R,R\Box R})\ci 
(R\Box \Sigma_{R\Box R,R\Box R}\Box R\Box R)\right)\\
\lb{j}&&{\bf j}=\left(R,\mu_R\ci(R\Box \mu_R\ci \Sigma_{R,R})\right).
\eea
The bimodule isomorphisms $l_{R^e},r_{R^e}$ and $a_{R^e}$ are constructed
as the unique morphisms in $\M$ for which
\bea 
\lb{lre}&&l_{R^e}\ci \sqcap(m,j\Box R^e)=(\mu_R\Box R)\ci \rho_m\ci
(R^{\Box 4}\Box \eta_R\Box R^{\Box 2})\\
\lb{rre}&&r_{R^e}\ci\sqcap(m,R^e\Box j)=(R\Box \mu_R)\ci \rho_m\ci 
(R^{\Box 5}\Box \eta_R\Box R)\\
\lb{are}&&a_{R^e}\ci\sqcap(m,m\Box R^e)=\sqcap(m,R^e\Box m)\ci (R\Box
R\Box \eta_R^{\Box 4}\Box R\Box R)\ci (R\Box R\Box \mu_R\Box R)\ci
\lambda_{m\Box R^e}\ci\nn
&&\qquad\qquad\qquad\qquad\quad(R\Box\eta_R\Box R^{\Box 7}).
\eea
The requirement that $\chi^*$ is a monoidal morphism from $({\bf
A},{\bf M},{\bf J})$ to $({\bf R^e},{\bf m},{\bf j})$ means the
existence of (bi-)module morphisms
\bea 
&&\sigma:{\bf m}\stac{{\bf R^e}\Box {\bf R^e}} (\chi^*\Box \chi^*)
\to \chi^*\stac{A} {\bf M}\nn
&&\iota:{\bf j}\to \chi^*\stac{A} {\bf J}
\nonumber\eea 
satisfying the conditions
\bea 
\lb{Lc} &&(\chi^*\stac{A} l_A)\ci [\sigma\stac{A\Box A} ({\bf J}\Box
{\bf A})]\ci [{\bf m}\stac {R^e\Box R^e} (\iota\Box \chi^*)]=
l_{R^e}\stac{R^e} \chi^*\\
\lb{Rc}&&(\chi^*\stac{A} r_A)\ci [\sigma\stac{A\Box A} ({\bf A}\Box
{\bf J})]\ci [{\bf m}\stac {R^e\Box R^e} (\chi^*\Box \iota)]=
r_{R^e}\stac{R^e} \chi^*\\
\lb{Cc}&&(\chi^*\stac{A} a_A)\ci [\sigma\stac{A\Box A} ({\bf M}\Box
{\bf A})]\ci [{\bf m}\stac {R^e\Box R^e} (\sigma\Box \chi^*)]=\nn
&&\qquad
[\sigma\stac{A\Box A} ({\bf A}\Box {\bf M})]\ci 
[{\bf m}\stac {R^e\Box R^e} (\chi^*\Box \sigma)]\ci
[a_{R^e}\stac{R^e\Box R^e\Box R^e} (\chi^*\Box\chi^*\Box\chi^*)].
\eea
The {\em strong} monoidality of $\chi^*$ means the additional
requirement that $\sigma$ and $\iota$ are isomorphisms.

\bigskip

Let us introduce the category ${\bf DS_R}(\M)$ as follows. The objects
are the pairs consisting of a pseudo-monoid $({\bf A,M,J})$ in $\BIM(\M)$
and a strong monoidal morphism from $({\bf A,M,J)}$ to the
canonical pseudo-monoid $({\bf R^e,m,j})$ which is induced by a 
monoidal morphism $\chi:{\bf R^e}\to {\bf A}$ in $\M$. Such pairs are
characterized by the collection of the data 
$({\bf A,M,J,R},\chi,\sigma,\iota)$.  

The morphisms  
$({\bf A},{\bf M},{\bf J},{\bf R},\chi,\sigma,\iota) \to  
({\bf A}^{\prime},{\bf M}^{\prime},{\bf J}^{\prime},{\bf R},
\chi^{\prime},$
$\sigma^{\prime},$ $\iota^{\prime})$ in ${\bf DS_R}(\M)$ are the morphisms
$\Phi:A\to A^{\prime}$ in $\M$ such that
\bea
\lb{PMmori} &&\Phi\ {\textrm{is\ a\ monoidal\ morphism}}\ {\bf A}\to
{\bf A}^{\prime} \ {\textrm{in}}\ \M\\
\lb{PMmorii} &&{\textrm{the }}\ {\bf A}-{\bf A}^{\prime}\ {\textrm{bimodule}}\
\Phi^*\     {\textrm{induced\ by}}\ \Phi \ {\textrm{is\ a\ monoidal\
morphism}}\nn 
&&\qquad ({\bf
    A^{\prime},M^{\prime},J^{\prime}})\to ({\bf A,M,J})\
{\textrm{in}}\ \BIM(\M) 
\\
\lb{PMmoriii} &&{\textrm{the\ composite\ of\ the\ monoidal\
    morphisms}}\ 
\Phi^*\ {\textrm{and}}\ \chi^*\ {\textrm{is\ equal\ to}}\  \chi^{\prime*}.
\eea
More explicitly the conditions (\ref{PMmorii}-\ref{PMmoriii}) mean that
$\Phi\ci\chi=\chi^{\prime}$ and
\bea\lb{tauform}&&\sigma^{\prime}\ci[\sigma^{-1}\stac{A\Box
    A}(\Phi^*\Box \Phi^*)]: 
M\stac{A\Box A} (\Phi^*\Box \Phi^*) \to \Phi^*\stac{A}
M^{\prime}\qquad {\textrm{and}}\\
\lb{jotaform}&&\iota^{\prime}\ci\iota^{-1}:J\to \Phi^*\stac{A} J^{\prime}
\eea
are (bi-) module isomorphisms satisfying the compatibility conditions
of the type (\ref{Lc}-\ref{Cc}).

\bigskip

From now on let ${\bf R}$ be a fixed monoid in $\M$, the $({\bf
R^e},{\bf m},{\bf j})$ the canonical pseudo-monoid (\ref{m}-\ref{j}).
For a monoidal morphism $\chi:{\bf R^e}\to {\bf A}$ let $s\colon
=\chi\ci(R\Box \eta_R)$ and $t\colon = \chi\ci(\eta_R\Box R)$
and denote the ${\bf R}$-${\bf R}$-bimodule (\ref{L}) by ${ L}$.

\bl Let ${\bf R}$ and ${\bf A}$ be monoids and $\chi:{\bf R^e}\to {\bf
A}$ a monoidal morphism in $\M$.
Then the left ${\bf R^e}$-modules  ${\bf m\stac{R^e\Box
R^e}}$ $ (\chi^*\Box \chi^*)$ and $\left(L\stac{R}
L,\lambda_{L\stac{R}L}\ci (R\Box \rho_{L\stac{R} L})\ci (R\Box
\Sigma_{R,L\stac{R}L})\right)$ in $\M$ are isomorphic. 
\el
\pr We construct the required isomorphism $\phi:m\stac{R^e\Box R^e}
(\chi^*\Box \chi^*)\to L\stac{R} L$ as the unique morphism which
satisfies
\be\lb{phi}
\phi\ci\sqcap(m,\chi^*\Box \chi^*)=\sqcap(L,L)\ci 
\lambda_{\chi^*\Box \chi^*}
%\mu_{A\Box A}\ci(\chi\Box \chi\Box A\Box A)
\ci (R\Box\eta_R\Box R\Box R\Box A\Box A).\ee
It is obviously a left ${\bf R^e}$-module  morphism.
It is an isomorphism with inverse the unique morphism ${\tilde \phi}$
for which
\be \lb{phiinv}
{\tilde \phi}\ci\sqcap(L,L)=\sqcap(m,\chi^*\Box\chi^*)\ci
(\eta_R^{\Box 3}\Box A\Box A).
\ee
\hfill $\Box$

For an object $({\bf A},{\bf M},{\bf J},{\bf R},\chi,\sigma,\iota)$ 
in ${\bf DS_R}(\M)$ let us use the isomorphism $\phi$ in (\ref{phi}) to
define the ${\bf A}$-${\bf A}\Box {\bf A}$-bimodule
\be \lb{ma}
{\bf M_A}=\left(L\stac{R}L, \lambda_{M_A}=\phi\ci \sigma^{-1}\ci
\lambda_M\ci(A\Box \sigma\ci \phi^{-1}), 
\rho_{M_A}=\phi\ci \sigma^{-1}\ci\rho_M\ci(\sigma\ci
\phi^{-1}\Box A\Box A)\right).\ee
It is easy to see that $\rho_{M_A}\ci[\sqcap(L,L)\Box A\Box
A]=\sqcap(L,L)\ci \mu_{A\Box A}$.
\bl \lb{xizeta} 
Let $({\bf A},{\bf M},{\bf J},{\bf R},\chi,\sigma,\iota)$ be an object in ${\bf
DS_R}(\M)$ and $M_A$ the ${\bf A}$-${\bf A}\Box {\bf A}$-bimodule
(\ref{ma}). Then the unique morphisms $\xi$ and $\zeta$ in $\M$ for
which the relations
\bea \lb{zeta}&&
\zeta\ci \sqcap(L,L\stac{R}L)=\sqcap(M_A,A\Box M_A)\ci
[\sqcap(L,L)\ci(\eta_A\Box\eta_A) \Box A\Box (L\stac{R}L)]\\
\lb{xi}&&
\xi\ci \sqcap(L\stac{R}L,L)=\sqcap(M_A,M_A\Box A)\ci
[\sqcap(L,L)\ci(\eta_A\Box\eta_A) \Box (L\stac{R}L)\Box A]
\eea
hold true are both isomorphisms and 
\be\lb{aform}
\zeta\ci \xi^{-1}=[\phi\ci \sigma^{-1}\stac{A\Box A}(A\Box \phi\ci
\sigma^{-1})] \ci a_A \ci
[\sigma\ci\phi^{-1}\stac{A\Box A}(\sigma\ci\phi^{-1}\Box A)] \ee
where $a_A$ is the coherent bimodule isomorphism satisfying
(\ref{Cc}). In particular, $\zeta\ci \xi^{-1}$ is an ${\bf A}$-${\bf
A}^{\Box 3}$-bimodule isomorphism.
\el
\pr The isomorphism property of $\zeta$ and $\xi$ is proven by
constructing their inverses as the unique morphisms ${\tilde\zeta}$
and ${\tilde\xi}$ in $\M$ for which
\bea \lb{zetainv}&&
{\tilde\zeta}\ci \sqcap(M_A,A\Box M_A)\ci[\sqcap(L,L)\Box A\Box M_A]=
\sqcap(L,L\stac{R}L)\ci \lambda_{A\Box M_A}\\
\lb{xiinv}&&
{\tilde\xi}\ci \sqcap(M_A,M_A\Box A)\ci[\sqcap(L,L)\Box M_A\Box A]=
\sqcap(L\stac{R}L,L)\ci \lambda_{M_A\Box A}.
\eea
The equation (\ref{aform}) follows from (\ref{zeta}) and 
\bea &&[\phi\ci\sigma^{-1}\stac{A\Box A} (A\Box \phi\ci\sigma^{-1})]
\ci a_A\ci 
[\sigma\ci\phi^{-1}\stac{A\Box A}(\sigma\ci\phi^{-1}\Box A)]\ci\xi\ci
\sqcap(L,L\stac{R} L)\ci (A\Box \sqcap(L,L))\nn
&&=[\phi\stac{A\Box A}(A\Box \phi)]\ci 
[a_{R^e}\stac{R^e\Box R^e\Box R^e} (\chi^*\Box \chi^*\Box \chi^*)]
\ci [\phi^{-1}\stac{A\Box A}(\phi^{-1}\Box A)]\ci \sqcap(M_A,M_A\Box
A)\nn
&&\qquad
\ci [\sqcap(L,L)\ci(\eta_A\Box \eta_A)\Box \sqcap(L,L)\Box A]\nn
&&=\sqcap(M_A,A\Box M_A)\ci [\sqcap(L,L)\ci(\eta_A\Box \eta_A)
\Box A\Box \sqcap(L,L)]
\nonumber\eea
which is the consequence of (\ref{xi}), (\ref{Cc}) and (\ref{are}).
\hfill $\Box$

For an object $({\bf A},{\bf M},{\bf J},{\bf R},\chi,\sigma,\iota)$ in ${\bf
DS_R}(\M)$ let us use the isomorphisms $\zeta$ and $\xi$ in
(\ref{zeta}-\ref{xi}) to define the ${\bf A}$-${\bf A}^{\Box
3}$-bimodule  
\bea \lb{m3}
&&{\bf M_3}=\left( L\stac{R} L\stac{R} L,
\lambda_{M_3}=\xi^{-1}\ci \lambda_{M_A\ot(M_A\Box A)}\ci (A\Box \xi)=
\zeta^{-1}\ci \lambda_{M_A\ot(A\Box M_A)}\ci (A\Box \zeta) \right. ,\nn
&&\qquad\quad
\left. 
\rho_{M_3}=\xi^{-1}\ci \rho_{M_A\ot(M_A\Box A)}\ci (\xi\Box A^{\Box 3})=
\zeta^{-1}\ci \lambda_{M_A\ot(A\Box M_A)}\ci (\zeta\Box A^{\Box 3})\right).
\eea
\bt The categories ${\bf Bgd_R}(\M)$ and ${\bf DS_R}(\M)$ are
equivalent. 
\et
\pr We start with the construction of a functor ${\cal
F}:{\bf DS_R}(\M) \to {\bf Bgd_R}(\M)$. It maps the object $({\bf
A},{\bf M},{\bf J},{\bf R},\chi,\sigma,\iota)$ of ${\bf DS_R}(\M)$ to the
bialgebroid 
$({\bf A},{\bf R},s,t,\gamma,\pi)$ where
\bea &&s=\chi\ci(R\Box \eta_R)\qquad t=\chi\ci(\eta_R\Box R)\\
&&\lb{gamma}\gamma=\lambda_{M_A}\ci[A\Box \sqcap(L,L)\ci(\eta_A\Box
\eta_A)]\\
&&\lb{pi}\pi=\iota^{-1}\ci\lambda_J\ci(A\Box \iota\ci\eta_R)
\eea
and the morphism $\lambda_{M_A}$ is introduced in (\ref{ma}).
Both $\gamma$ and $\pi$ are ${\bf R}$-${\bf R}$-bimodule morphisms  
by the left ${\bf R^e}$-module morphism property of $\iota$ and of 
$\phi\ci \sigma^{-1}$,
respectively. The axioms (\ref{bgdi}-\ref{bgdv}) follow easily from
the definitions (\ref{gamma}-\ref{pi}) using the identities $\pi\ci
s=R$ and 
\be \lb{lrac} \rho_{M_A}\ci (\gamma\Box A\Box A)=\lambda_{M_A}\ci
(A\Box \sqcap(L,L)).\ee
The non-trivial properties are the coassociativity of $\gamma$ and the
counit property of $\pi$. The coassociativity follows from the Lemma
\ref{xizeta} as follows. Use the first forms of $\lambda_{M_3}$ and of
$\rho_{M_3}$ in (\ref{m3}) and the identity (\ref{lrac}) to conclude
that
$$\lambda_{M_3}\ci[A\Box \sqcap(L\stac{R}L,L)\ci(\sqcap(L,L)\Box A)]
=\rho_{M_3}\ci [(\gamma\stac{R} L)\ci \gamma\Box A\Box A\Box A]$$
and the second forms  of $\lambda_{M_3}$ and of $\rho_{M_3}$ in
(\ref{m3}) and (\ref{lrac}) to show that 
$$\lambda_{M_3}\ci[A\Box \sqcap(L,L\stac{R}L)\ci(A\Box \sqcap(L,L))]
=\rho_{M_3}\ci [(L\stac{R} \gamma)\ci \gamma\Box A\Box A\Box A].$$
This results $(L\stac{R} \gamma)\ci \gamma=\lambda_{M_3}\ci
\{A\Box \sqcap(L\stac{R}L,L)\ci[\sqcap(L,L)\ci(\eta_A\Box \eta_A)\Box
\eta_A]\} =(\gamma\stac{R} L)\ci \gamma$.

The counit property is proven by a similar argument. Introduce the
left ${\bf A}$-module 
\be\lb{ja} {\bf J_A}=\left(R,\iota^{-1}\ci\lambda_J\ci (A\Box \iota)\equiv
\pi\ci\mu_A\ci(A\Box s)\right).\ee
One checks that the unique morphisms in $\M$ which satisfy
\bea \lb{kappa1}&&
\kappa_1\ci\sqcap(M_A,J_A\Box A)=(\pi\stac{R} L)\ci \rho_{M_A}\ci
[(L\stac{R}L)\Box s\Box A]\quad {\textrm{and}}\\
\lb{kappa2}&&
\kappa_2\ci\sqcap(M_A,A\Box J_A)=(L\stac{R} \pi)\ci \rho_{M_A}\ci
[(L\stac{R}L)\Box A\Box t]\eea
are $\kappa_1=l_A\ci [\sigma\ci\phi^{-1}\stac{A\Box
A}(\iota\Box A)]$ and  $\kappa_2=r_A\ci[\sigma\ci\phi^{-1}\stac{A\Box
A}(A\Box \iota)]$, respectively, where $l_A$ and $r_A$ are the coherent bimodule
isomorphisms satisfying the conditions (\ref{Lc}) and (\ref{Rc}). Now
using the left ${\bf A}$-module morphism property of $\kappa_1$ and
(\ref{lrac}) we obtain
\bea 
(\pi\stac{R}L)\ci \gamma&=&(\pi\stac{R}L)\ci \rho_{M_A}\ci
[\gamma\Box\mu_{A\Box A}\ci(\eta_A\Box\eta_A\Box s\ci\eta_R\Box
\eta_A)]=\nn
&&\kappa_1\ci \lambda_{M_A\ot(J_A\Box A)}\ci
\{A\Box \sqcap(M_A,J_A\Box A)\ci[\sqcap(L,L)\ci(\eta_A\Box\eta_A)\Box
\eta_R\Box \eta_A]\}=\nn
&&\mu_A\ci[A\Box (\pi\stac{R}L)\ci\sqcap(L,L)\ci(\eta_A\Box
\eta_A)]=A\eea
and from the left ${\bf A}$-module morphism property of $\kappa_2$ it
follows that $(L\stac{R}\gamma)\ci \gamma=A$. This finishes the
construction of the functor ${\cal F}$ on the objects.

The functor ${\cal F}$ acts on the morphisms as the identity map. For
a morphism $\Phi:({\bf A},{\bf M},{\bf J},{\bf R},\chi,\sigma,$
$\iota)  \to({\bf
 A}^{\prime}, {\bf M}^{\prime},{\bf J}^{\prime},{\bf R},\chi^{\prime},
\sigma^{\prime},\iota^{\prime})$  in ${\bf DS_R}(\M)$ the conditions
(\ref{bgdmori}-\ref{bgdmorii}) follow from the fact that $\Phi$ is a
bimodule morphism $L\to L^{\prime}$ and a left ${\bf R^e}$-module
morphism $\chi^*\to \chi^{\prime
*}$ and the left ${\bf A }$-module morphism property of the morphisms
(\ref{tauform}-\ref{jotaform}).

Let us turn to the construction of the functor ${\cal G}:
{\bf Bgd_R}(\M)\to {\bf DS_R}(\M)$. It maps the object ${\cal A}=({\bf
A},{\bf R},s,t,\gamma,\pi)$ in ${\bf Bgd_R}(\M)$ to $({\bf A},{\bf
M_A},{\bf J_A},{\bf R},\chi, \phi,R)$ where 
${\bf M_A}$ is the ${\bf A}$-${\bf A}\Box {\bf A}$-bimodule defined 
on the object $L\stac{R} L$ as in (\ref{rho}-\ref{lambda}), the ${\bf
J_A}$ is the left ${\bf A}$-module
$$\left(R,\pi\ci\mu_A\ci(A\Box R)\right)$$
and $\chi\colon =\mu_A\ci (s\Box t)$ and the morphism $\phi$
has been introduced in (\ref{phi}).

We claim that $({\bf A},{\bf M_A},{\bf J_A})$ is a pseudo-monoid in 
$\BIM(\M)$ by constructing the coherent bimodule isomorphisms
$l_A,r_A$ and $a_A$. The $l_A$ and $r_A$ are the unique morphisms in
$\M$ which satisfy
\bea \lb{la}&&
l_A\ci \sqcap(M_A,J_A\Box A)=(\pi\stac{R} L)\ci \rho_{M_A}\ci [(L\stac{R}L)\Box
s\Box A]\qquad {\textrm{and}}\\
\lb{ra}&& 
r_A\ci \sqcap(M_A,A\Box J_A)=(L\stac{R} \pi)\ci \rho_{M_A}\ci [(L\stac{R}L)\Box
A\Box t],\eea
respectively. Both are isomorphisms with inverses
\bea &&
l_A^{-1}=\sqcap(M_A,J_A\Box A)\ci[\sqcap(L,L)\ci(\eta_A\Box
\eta_A)\Box \eta_R\Box A]\nn
&&r_A^{-1}=\sqcap(M_A,A\Box J_A)\ci[\sqcap(L,L)\ci(\eta_A\Box \eta_A)
\Box A\Box \eta_R].
\nonumber\eea
Their right ${\bf A}$-module morphism property is
obvious. The left ${\bf A}$-module morphism property follows 
from (\ref{bgdv}), the left ${\bf R}$-module morphism property of
$\gamma$, the axioms (\ref{bgdii}) and (\ref{bgdiii}) and the  counit
property of $\pi$. 

In order to construct $a_A$ introduce the morphisms $\zeta$ and $\xi$
with the same formulae (\ref{zeta}-\ref{xi}) as in Lemma \ref{xizeta}
with the only difference that $M_A$ stands now for the ${\bf A}$-${\bf
A}\Box{\bf A}$-bimodule (\ref{rho}-\ref{lambda}). Just the same way as
in the Lemma they are shown to be isomorphisms. Set 
\be \lb{aa} a_A\colon = \zeta\ci\xi^{-1}.\ee
Its bimodule morphism property follows from the bimodule morphism
property of
\bea &&\zeta:(L\stac{R} L\stac{R} L,\lambda_3,\rho_3)\to 
{\bf M_A\stac{A\Box A} (M_A\Box A)}\qquad {\textrm{and}}\nn
&&\xi: (L\stac{R} L\stac{R} L,\lambda_3^{\prime},\rho_3)\to 
{\bf M_A\stac{A\Box A} (A\Box M_A)}\nonumber \eea
(where the morphisms $\lambda_3,\lambda_3^{\prime}$ and $\rho_3$ are
the unique morphisms in $\M$ for which
\bea &&
\rho_3\ci[\sqcap(L\stac{R}L,L)\Box A^{\Box 3}]=
\sqcap(L\stac{R}L,L)\ci \rho_{M_A\Box A}\quad \Leftrightarrow\\
&&\rho_3\ci[\sqcap(L,L\stac{R}L)\Box A^{\Box 3}]=
\sqcap(L,L\stac{R}L)\ci \rho_{A\Box M_A}\\
\lb{lambda32}&&
\lambda_3\ci[A\Box \sqcap(L,L\stac{R}L)\ci(A\Box \sqcap(L,L))]=
\rho_3\ci [(L\stac{R} \gamma)\ci \gamma\Box A^{\Box 3}]\\
&&\lambda_3^{\prime}\ci[A\Box \sqcap(L\stac{R}L,L)\ci(\ \sqcap(L,L)\Box
A)\ ]=
\rho_3\ci [(\gamma\stac{R} L)\ci \gamma\Box A^{\Box 3}]\qquad )
\lb{lambda31}
\eea
and the coassociativity of $\gamma$.

The proof of the pseudo-monoid property of $({\bf A},{\bf M_A},{\bf
J_A})$ is completed by the lengthy but
straightforward check of the triangle and pentagon conditions. 

The fact that the bimodule isomorphisms
\bea &&\phi: {\bf m}\stac{R^e\Box R^e} (\chi^*\Box \chi^*)\to
\chi^*\stac{A} {\bf M_A} \nn
&& R:{\bf j}\to \chi^*\stac{A} {\bf J_A}
\nonumber\eea
satisfy the conditions (\ref{Lc}-\ref{Cc}) follows from
the identities
(\ref{la}) and (\ref{lre}), (\ref{ra}) and (\ref{rre}), (\ref{aa}) and
(\ref{are}) -- after some playing with the coequalizers. 

The functor ${\cal G}$ acts on the morphisms as the identity map. One
checks that
for a morphism $\Phi:({\bf A},{\bf R},s,t,\gamma,\pi)$ $\to ({\bf
A}^{\prime},{\bf R},s^{\prime},t^{\prime},\gamma^{\prime},
\pi^{\prime})$ in ${\bf Bgd_R}(\M)$ the conditions
(\ref{PMmori}-\ref{PMmoriii}) follow from
(\ref{bgdmori}-\ref{bgdmorii}). 

The composite functor ${\cal F}{\cal G}$ is the identity functor ${\bf
Bgd_R}(\M)$, while ${\cal G}{\cal F}$ is naturally equivalent to the
identity functor ${\bf DS_R}(\M)$. The image of the object $({\bf
A},{\bf M},{\bf J},{\bf R},\chi,\sigma,\iota)$ in  ${\bf DS_R}(\M)$ under
${\cal G}{\cal F}$ is the object $({\bf A},{\bf M_A},{\bf J_A},{\bf
R},\chi,\phi,R)$ where ${\bf M_A}$ is the ${\bf A}-{\bf A}\Box {\bf
  A}$-bimodule defined in
(\ref{ma}), ${\bf J_A}$ is the left ${\bf A}$-module (\ref{ja}) and
$\phi$ is the morphism (\ref{phi}). The required natural equivalence at the
object $({\bf A},{\bf M},{\bf J},{\bf R},\chi,\sigma,\iota)$ is the 
identity morphism $A$. As a matter of fact the ${\bf A}$-${\bf
A}$-bimodule $(A,\mu_A,\mu_A)$ with the coherent bimodule isomorphisms
$\sigma\ci \phi^{-1}:{\bf M_A}\to {\bf M}$ and
$\iota:{\bf J_A}\to {\bf J} $ is a strong monoidal morphism $({\bf
A,M,J})\to ({\bf A,M_A,J_A})$ since the
coherence conditions (\ref{Lc}) and (\ref{Rc}) hold true by the
observation that the unique morphisms satisfying
(\ref{kappa1}-\ref{kappa2}) are $l_A\ci[\sigma\ci\phi^{-1}\stac{A\Box
A} (\iota\Box A)]$ and $r_A\ci[\sigma\ci\phi^{-1}\stac{A\Box
A} (A\Box \iota)]$, respectively, and the definitions
(\ref{la}-\ref{ra}). The condition (\ref{Cc}) follows from Lemma
\ref{xizeta} and (\ref{aa}). 
\hfill $\Box$

\bigskip

The definition (\ref{bgdmori}-\ref{bgdmorii}) of bialgebroid morphisms
is very restrictive. In the case of bialgebroids in $\M_k$ more
general morphisms have been introduced in \cite{Schau2}
and also in \cite{Szlach}. The idea of \cite{Szlach} can be applied
also in our context. One can introduce a bigger category ${\bf
DS}(\M)$ the objects of which are the objects of all the categories
${\bf DS_R}(\M)$ as ${\bf R}$ runs through the monoids in $\M$. The
definition of the morphisms follows Street's definition of the
1-cells in the bicategory of monads \cite{Street}. Then
the object functions of the functors ${\cal F}$ and ${\cal G}$
constructed in the proof can be used to define a bigger category ${\bf
Bgd}(\M)$ of bialgebroids in $\M$. As a guiding principle they are
required to be the object functions of equivalence functors ${\cal
F}^{\prime}: {\bf Bgd}(\M) \to {\bf DS}(\M)$ and ${\cal
G}^{\prime}: {\bf DS}(\M)\to {\bf Bgd}(\M)$, respectively.
We do not consider this question in more detail here.

\bigskip

The axioms of the internal left bialgebroid in the Definition \ref{bgd} are not
invariant under the change of the monoid ${\bf A}$ to the opposite ${\bf
A}^{op}$. As the opposite version of the left bialgebroid the right bialgebroid
is introduced as follows.
\bd \lb{rbgd}
A {\em right bialgebroid} in the symmetric monoidal category $(\M,\Box, U)$
with coequalizers
consists of the data $({\bf A},{\bf R},s,t,\gamma,\pi)$ where ${\bf
A}=(A,\mu_A,\eta_A)$ and ${\bf R}=(R,\mu_R,\eta_R)$ are monoids in $\M$, the
$s:{\bf R}\to {\bf A}$ and the $t:{\bf R^{op}}\to {\bf A}$ are monoidal
morphisms such that
$$\mu_A\ci(s\Box t)=\mu_A\ci\Sigma_{A,A}\ci(s\Box t)$$
and ${\bf K}=(K,\gamma,\pi)$ is a comonoid in ${_{\bf R}\M_{\bf R}}$ where $K$
is the ${\bf R}$-${\bf R}$-bimodule 
$$ \left(A,\mu_A\ci\Sigma_{A,A}\ci(t\Box A),\mu_A\ci (A\Box s)\right). $$
They are subject to the axioms
\bea &&\lambda_{N_A}\ci (\eta_A\Box t\Box \gamma)=\lambda_{N_A}\ci (s\Box 
\eta_A\Box \gamma) \nn
&&\gamma\ci \eta_A=\sqcap(K,K)\ci (\eta_A\Box \eta_A)\nn
&&\gamma\ci\mu_A=\rho_{N_A}\ci (\gamma\Box A)\nn
&&\pi\ci\eta_A=\eta_R\nn
&&\pi\ci\mu_A\ci(s\ci\pi\Box A)=\pi\ci\mu_A=\pi\ci\mu_A\ci(t\ci\pi\Box A)
\eea
where we introduced the ${\bf A\Box A}$-${\bf A}$-bimodule ${\bf N_A}=
(K\stac{R} K, \lambda_{N_A},\rho_{N_A})$ with the unique morphisms in $\M$ for
which
\bea &&\lambda_{N_A}\ci(A\Box A\Box \sqcap(K,K))=\sqcap(K,K)\ci
\mu_{A\Box A}\nn
&&\rho_{N_A}\ci(\sqcap(K,K)\Box A)=\lambda_{N_A}\ci(A\Box A\Box \gamma).
\nonumber\eea
\ed

\section{The formulation of internal bialgebroids in terms of
entwining structures and corings}
\lb{relate}
\setc{0}

The axioms of the bialgebra over a commutative ring $k$ have been
related to entwining structures in (\cite{BreWis}, Proposition 33.1)
and to corings in (\cite{Wis}, Proposition 5.2). Our aim in this 
section is to give similar relation between bialgebroids, entwining
structures and corings in symmetric monoidal categories with
coequalizers. 

Throughout the section $(\M,\Box,U,\Sigma)$ is a symmetric monoidal
category with coequalizers, ${\bf R}=(R,\mu_R,\eta_R)$ and ${\bf
A}=(A,\mu_A,\eta_A)$ are monoids in $\M$ and $\chi:{\bf R^e}\to {\bf
A}$ is a monoidal morphism. Let us introduce the ${\bf R}$-${\bf
R}$-bimodules
\bea &&L=\left(A,\lambda_L=\mu_A\ci(s\Box
A),\rho_L=\mu_A\ci\Sigma_{A,A}\ci(A\Box t)\right)\nn
&&S=\left(A,\lambda_S=\mu_A\ci(s\Box A),\rho_S=\mu_A\ci (A\Box s)
\right)
\nonumber\eea  
where $s=\chi\ci(R\Box \eta_R)$ and $t=\chi\ci(\eta_R\Box
R)$. (Notice that the $\M$-objects $L\stac{R} L$ and $L\stac{R} S$
coincide.) 
Introducing the unique morphism $\mu_S:S\stac{R} S\to S$ for
which 
$$\mu_S\ci \sqcap(S,S)=\mu_A $$
we have a monoid ${\bf S}=(S,\mu_S,s)$ in ${_{\bf R}\M_{\bf R}}$, the
monoidal category of ${\bf R}$-${\bf R}$-bimodules in $\M$.

Let ${\bf L}=(L,\gamma,\pi)$ be a comonoid in ${_{\bf R}\M_{\bf R}}$ 
and $\rho_{M_A}: (L\stac{R} L)\Box A \Box A\to  (L\stac{R}
L)$ the morphism (\ref{rho}) in $\M$.
\bt \lb{entw-bgd}
The following are equivalent:
\begin{trivlist}
\item {\em i)} The comonoid ${\bf L}$ obeys the properties
\bea \lb{1i} &&\pi\ {\textrm{is\ an\ epimorphism}}\\
\lb{1ii}&&\gamma\ci\rho_S=\rho_{M_A}\ci(\gamma\Box s\Box \eta_A)\\
&&{\textrm{the\ monoid\ }} {\bf S}, {\textrm{the\ comonoid\
}} {\bf L}\ {\textrm{and\ the\ unique\ morphism\ }} \psi: {\bf
S}\stac{R} {\bf L}\to {\bf L}\stac{R}{\bf S}
%\nn&&
{\textrm{\ for \ which}}\nn
&&\qquad\qquad \psi\ci\sqcap(S,L)=\rho_{M_A}\ci(\gamma\Box A\Box
\eta_A) \nn
\lb{1iii}&&{\textrm{form\ a\ left\ entwining\ structure\ in\ }} {_{\bf
R}\M_{\bf R}}.
\eea
\item {\em ii)} The comonoid ${\bf L}$ satisfies the properties
(\ref{1i}), (\ref{bgdi}) and
\be \lb{cor} 
{\bf C}=\left(L\stac{R} S, \lambda_{M_A},\rho_{M_A}\ci [(L\stac{R} S)\Box
\eta_A\Box A], \gamma\stac{R} S,\pi\stac{R} S\right)\ee
is an ${\bf A}$-coring in $\M$, where $\lambda_{M_A}:A\Box (L\stac{R}
S)\to (L\stac{R} S)$ is the morphism (\ref{lambda}) in $\M$.
\item {\em iii)} $({\bf A},{\bf R},s,t,\gamma,\pi)$ is a left bialgebroid
in $\M$.
\end{trivlist}
\et
\br The implication $iii)\Rightarrow ii)$ of Theorem \ref{entw-bgd}
follows from (\cite{BreCaenMil}, 
Proposition 4.1) in the case when $\M$ is the category of modules over
a commutative ring.
\er
\br \lb{rem} The $i\Leftrightarrow iii)$ part of the Theorem
\ref{entw-bgd} is a generalization of 
(\cite{BreWis}, Proposition 33.1 (1)) where no analogues of the conditions
(\ref{1i}) and (\ref{1ii}) appear. 

We need the condition (\ref{1ii}) in order for the definition
(\ref{1iii}) of the morphism $\psi$ to make sense. In the context of
\cite{BreWis} one works in the category of modules over a commutative
ring $k$. In this case the condition (\ref{1ii}) is trivial in the
sense that it is identical to the $k$-module morphism property of
$\gamma$. 

The role of (\ref{1i}) is different. The starting point in
\cite{BreWis} is an algebra
$(A,\mu,\eta)$ and a coalgebra $(A,\gamma,\pi)$ over a commutative
ring $k$ which are entwined by the $k$-linear map
$$\psi:A\stac{k}A\to A\stac{k}A\qquad a\stac{k} b\mapsto a_{(1)}
b\stac{k} a_{(2)} $$
-- where we used the Sweedler's convention for the coproduct. It is
proven in \cite{BreWis} that both $\pi$ and $\gamma$ are
multiplicative and $\gamma$ is unital. The reader may check that the
unitality of $\pi$ is equivalent, however, to the faithfulness of $A$
as a $k$-module, what is further equivalent to the injectivity of
$\eta$ and also equivalent to the surjectivity of $\pi$.

The situation in our more general context is similar. All bialgebroid
axioms (\ref{bgdi}-\ref{bgdv}) follow from the assumptions
(\ref{1ii}-\ref{1iii}) except (\ref{bgdiv}). As a matter of fact under
the assumptions (\ref{1ii}-\ref{1iii}) the following are equivalent:
\bea && \pi\ {\textrm{is\ an\ epimorphism}}\nn
     \Leftrightarrow\quad && s\ {\textrm{is\ a\ monomorphism}}\nn
     \Leftrightarrow\quad && \pi\ci\eta_A=\eta_R.
\nonumber\eea
\er
{\em Proof of Theorem \ref{entw-bgd}.} 

$i)\Rightarrow iii)$ By the assumption (\ref{1ii})
the definition (\ref{1iii}) of $\psi$ makes sense. Since $\psi$ is
assumed to be a right ${\bf R}$-module morphism the axiom (\ref{bgdi})
is satisfied. From the assumption that $({\bf S},{\bf L},\psi)$
satisfies (\ref{entwi}) it follows that the axiom (\ref{bgdii}) holds
true, by (\ref{entwiii}) so does (\ref{bgdiii}). The axiom
(\ref{bgdiv}) follows from the assumption (\ref{1i}) while
(\ref{bgdv}) is the consequence of (\ref{entwii}) and the identity 
$$ \pi\ci (\pi\stac{R} S)\ci \rho_{M_A}\ci[ (L\stac{R} L)\Box A\Box
\eta_A]= \pi\ci \mu_A\ci [(L\stac{R} \pi)\Box A] $$
which follows from the ${\bf R}$-${\bf R}$-bimodule morphism property
of $\pi$.

$iii)\Rightarrow i)$ It follows from (\ref{bgdiv}) and the left ${\bf
R}$-module morphism property of $\pi$ that $\pi\ci s=R$ hence the
condition (\ref{1i}) holds true. The condition (\ref{1ii}) follows
from (\ref{bgdiii}), (\ref{bgdii}) and the left ${\bf R}$-module
morphism property of $\gamma$. The left ${\bf R}$-module morphism
property of $\psi$ follows from the left ${\bf R}$-module morphism property of
$\gamma$ and its right ${\bf R}$-module morphism property follows from
(\ref{bgdi}). The axiom (\ref{entwi}) follows from (\ref{bgdii}). The
(\ref{bgdv}) implies that
\be\lb{auxi22}(\pi\stac{R} S)\ci\rho_{M_A}\ci[(L\stac{R}L)\Box A\Box\eta_A]=
(\pi\stac{R} S)\ci\rho_{M_A}\ci[(L\stac{R}L)\Box s\ci\pi\Box\eta_A].\ee
Applying (\ref{auxi22}), the assumption (\ref{1ii}) and the counit
property of $\pi$ one proves (\ref{entwii}).

In order to prove (\ref{entwiii}) use the identity
$$
(L\stac{R}\mu_S)\ci (\psi\stac{R} S)\ci \sqcap(S,L\stac{R} S)\ci
(A\Box \rho_{M_A})=\rho_{M_A}\ci (\lambda_{M_A}\Box A\Box A)$$
together with (\ref{bgdiii}). 

The axiom (\ref{bgdiii}) implies that both $\gamma\stac{R} L$ and
$L\stac{R}\gamma$ are left ${\bf A}$-module morphisms $(L\stac{R}
L,\lambda_{M_A})\to (L\stac{R} L\stac{R} L,\lambda_3)$ where $\lambda_3$
has been introduced in (\ref{lambda31}). Furthermore
$$ (L\stac{R}\psi)\ci (\psi\stac{R}L)\ci \sqcap(S,L\stac{R}
L)=\lambda_3\ci \{A\Box \sqcap(L\stac{R} L,L)\ci[(L\stac{R} L)\Box
\eta_A]\}.$$
These observations imply (\ref{entwiv}).

$i)\Rightarrow ii)$ The axiom (\ref{bgdi}) holds true by the right
${\bf R}$-module morphism property of $\psi$.

Since $\M$ is a monoidal category with coequalizers so is ${_{\bf
R}\M_{\bf R}}$. 
Using the method of the Example \ref{entw-coring} we can construct an
${\bf S}$-coring in  ${_{\bf R}\M_{\bf R}}$ as
\be\lb{Scor}({\bf L}\stac{R} {\bf S},\Lambda,L\stac{R}\mu_S,\gamma\stac{R}
S,\pi\stac{R} S)\ee
where $\Lambda=(L\stac{R} \mu_S)\ci (\psi\stac{R} S)$.

The forgetful functor ${G^R}:{_{\bf R}\M_{\bf R}}\to \M$ is a
monoidal functor preserving the module tensor product hence by Proposition
\ref{forget} it maps the ${\bf S}$-coring (\ref{Scor}) in ${_{\bf
R}\M_{\bf R}}$ to an ${\bf A}$-coring in $\M$. It is straightforward to
check that it is the ${\bf A}$-coring (\ref{cor}).

$ii)\Rightarrow i)$ The unitality of the bimodule
(\ref{rho}-\ref{lambda}) implies that $\gamma\ci\eta_A=\sqcap(L,L)\ci
(\eta_A\Box\eta_A)$ and by the left ${\bf R}$-module morphism property
of $\gamma$ also $\gamma\ci s=\sqcap(L,L)\ci (s\Box\eta_A)$. Combining
this result with the assumption that $(L\stac{R} S,\lambda_{M_A})$ is
a left ${\bf A}$-module in $\M$, one proves the condition (\ref{1ii}).

This makes it sensible to introduce the unique morphism
$\Lambda:S\stac{R} L\stac{R} S \to 
L\stac{R} S$ in $\M$ with the requirement that
\be\lb{Lambda} \Lambda\ci \sqcap(S,L\stac{R} S)=\lambda_{M_A}.\ee
Observe that the ${\bf A}$-coring (\ref{cor}) in $\M$ is the image of
the ${\bf S}$-coring
$$ ({\bf L}\stac{R} {\bf S},\Lambda,L\stac{R}\mu_S,\gamma\stac{R}
S,\pi\stac{R} S)$$
in ${_{\bf R}\M_{\bf R}}$ under the forgetful functor ${G^R}:{_{\bf
R}\M_{\bf R}}\to \M$. Indeed, since ${G^R}$ is a monoidal functor preserving
the module tensor product and it is faithful and $G^R_2=\sqcap(\ ,\ )$ is epi,
the functor $G^R$ reflects corings by Proposition \ref{forget}. 

The morphism $\psi\colon = \Lambda\ci (S\stac{R} L\stac R
s): S\stac{R} L\to L\stac{R} S$ in $\M$ satisfies
$$\Lambda=(L\stac{R} \mu_S)\ci (\psi\stac{R} S) $$
hence we can apply the ``only if'' part of the statement in the Example
\ref{entw-coring} to conclude the claim.
\hfill $\Box$

It is straightforward to prove the `right handed' version of the 
Theorem \ref{entw-bgd} that is to formulate the axioms of the right
bialgebroid in $\M$ in terms of a  coring
in $\M$ and alternatively, in terms of a right entwining structure 
in ${_{\bf R}\M_{\bf R}}$.   

\bigskip

The construction of associating a coring (\ref{cor}) to a bialgebroid can be
applied to a bialgebra $H$ over a commutative ring $k$. As it is known
from (\cite{Nak}, Theorem 1.1) the resulting coring is Galois if and
only if $H$ is a Hopf algebra. Motivated by this result in the rest of
the section we investigate under what conditions on the bialgebroid
$({\bf A},{\bf R},s,t,\gamma,\pi)$ is the ${\bf A}$-coring (\ref{cor})
Galois. 

\begin{figure}[h]
\psfrag{3}{\Large$R$}
\psfrag{1}{\Large$A$}
\psfrag{2}{\Large$L\stac{R} S$}
\psfrag{6}{$t$}
\psfrag{4}{$\rho_{M_A}\ci [\sqcap(L,S)\ci(\eta_A\Box \eta_A)\Box \eta_A\Box A]$}
\psfrag{5}{$\lambda_{M_A}\ci [A\Box \sqcap(L,S)\ci(\eta_A\Box \eta_A)]$}
\begin{center}
{\resizebox*{12.5cm}{!}{\includegraphics{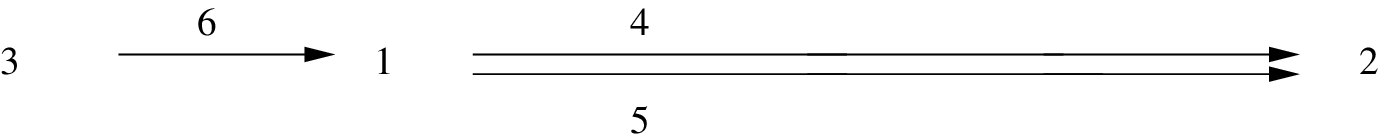}}
}
\end{center}
\caption{The coinvariant subobject in the coring (\ref{cor})}
\label{eqbgd}
\end{figure}

The Galois property of corings is defined w.r.t. a given group-like
morphism. Let us fix the group-like morphism $\sqcap(L,S)\ci(\eta_A\Box
\eta_A)$ for the coring (\ref{cor}). The corresponding coinvariant
subobject of $A$ is identified by the following
\bl \lb{coinv}
The diagram on the Figure \ref{eqbgd} is an equalizer diagram in $\M$.
\el
\pr The commutativity of the diagram on Figure \ref{eqbgd} follows
from the observation that $\rho_{M_A}\ci [\sqcap(L,S)\ci(\eta_A\Box
\eta_A)\Box \eta_A\Box A]=\sqcap(L,S)\ci (\eta_A\Box 
A)$ and $\lambda_{M_A}\ci [A\Box \sqcap(L,S)\ci(\eta_A\Box
\eta_A)]=\gamma$, the axiom (\ref{bgdii}) and the right ${\bf R}$-module
morphism property of $\gamma$.

The universality of the morphism $t:R\to A$ follows since for any
object $X$ and morphism $f:X\to A$ in $\M$ such that
$\sqcap(L,S)\ci(\eta_A\Box f)=\gamma\ci f$ the
morphism $f$ factorizes through $R$ in the unique way $f=t\ci\pi\ci
f$.
\hfill $\Box$

Recall that the Galois property of an internal coring is equivalent to
the isomorphism property of the morphism $\kappa$ defined by
(\ref{kappaform}). Substituting the group-like morphism
$\sqcap(L,S)\ci (\eta_A\Box \eta_A)$ and the equalizer in the Lemma
\ref{coinv} the definition (\ref{kappaform}) takes the form
\be \lb{bgdkappa} 
\kappa\ci \sqcap(t_*,t^*)=\rho_{M_A}\ci(\gamma\Box \eta_A\Box A).\ee
The morphism $\kappa$ introduced by (\ref{bgdkappa}) has been studied
by P. Schauenburg \cite{Schau} in the case of bialgebroids $({\bf
A},{\bf R},s,t,\gamma,\pi)$ in the category $\M_k$. In that case the
isomorphism property of 
$\kappa$ was shown to be equivalent to the strong right closed
property of the forgetful functor $_{\bf A}\M\to {_{\bf R}\M_{\bf
R}}$ hence it was proposed as the definition of the $\times_R$-Hopf
algebra. 
\bc Let ${\cal A}=({\bf A},{\bf R},s,t,\gamma,\pi)$ be a left bialgebroid 
in $\M$. The corresponding ${\bf A}$-coring (\ref{cor}) is Galois 
(with respect to the coinvariant subobject $R$ of $A$ corresponding to
the group-like morphism $\gamma\ci \eta_A$) if and
only if the morphism $\kappa$ in $\M$ introduced in (\ref{bgdkappa})
is an isomorphism. In this case we may call ${\cal A}$ a {\em
$\times_R$-Hopf algebra}.
\ec

\section{Duality of internal bialgebroids}

The entwining structure is a self-dual notion in the sense of Remark
\ref{entwdual}.
Therefore the the description of the internal bialgebroid in the
part $i)$ of Theorem \ref{entw-bgd} provides a natural framework 
for the study of its duality.
That is one can ask whether the dual of the
entwining structure in the part $i)$ of Theorem \ref{entw-bgd} -- if exists 
-- satisfies analogous properties to (\ref{1i}-\ref{1iii}) hence determines 
a dual bialgebroid. 
The answer is given by the following
\bt Let $(\M,\Box,U,\Sigma)$ be a symmetric monoidal category with
coequalizers. Let ${\cal A}=({\bf A},{\bf R},s,t,\gamma,\pi)$ be a left
bialgebroid in $\M$ and $({\bf S},{\bf L},\psi)$ the entwining structure in
${_{\bf R}\M_{\bf R}}$ as in the part $i)$ of Proposition
\ref{entw-bgd}. Suppose that both ${\bf R}$-${\bf R}$ bimodules $S$ and $L$
possess right duals in ${_{\bf R}\M_{\bf R}}$. Then the dual bimodules have
the form
\bea \lb{Lr} 
&&L^r=\left(B,\mu_B\ci(\pi^r\Box B),\mu_B\ci (B\Box \pi^r)\right)\\
\lb{Sr} &&S^r=\left(B, \mu_B\ci \Sigma_{B,B}\ci(q\Box B),\mu_B\ci(B\Box \pi^r)
\right),
\eea
respectively, where $B$ is some object and $q$ is some morphism in
$\M$ 
and  ${\bf B}$ is the monoid $(B,\gamma^r\ci \sqcap(L^r,L^r),\pi^r\ci
\eta_R)$ in $\M$. The right entwining structure $({\bf L}^r,{\bf S}^r,\psi^r)$
obeys the properties
\bea && s^r\ {\textrm{is\ an\ epimorphism}}\\
\lb{psir}&& \psi^r\ci \sqcap(S^r,L^r)=\lambda_{N_A}\ci (\eta_B\Box B
\Box \mu^r)
\eea
where $\lambda_{N_A}$ is the unique morphism in $\M$ for which $\lambda_{N_A}
\ci (B\Box B\Box \sqcap(L^r,S^r))=\sqcap(L^r,S^r)\ci\mu_{B\Box
B}$. Equivalently, $({\bf B},{\bf R}, \pi^r,q,\mu^r,s^r)$ is a right
bialgebroid in $\M$.
\et
\pr Let us denote the evaluation and coevaluation morphisms by
\bea &\cup_S:S\stac{R} S^r\to R\qquad &\cap_S:R\to S^r\stac{R} S\nn
&\cup_L :L\stac{R} L^r\to R\qquad &\cap_L:R\to L^r\stac{R} L.
\nonumber\eea
It is obvious that $L^r$ is of the form (\ref{Lr}) since ${\bf
L}^r=(L^r,\gamma^r,\pi^r)$ is a monoid in $\M$. Since $\lambda_L=\lambda_S$,
the $\M$-object parts $B$ of $L^r$ and $B^{\prime}$ of $S^r$ are isomorphic
via
\bea &&t\ \ \ \! \colon =(S^r\stac{R} \cup_L)\ci \sqcap(S^r\stac{R} L,L^r)\ci
(\cap_S\ci \eta_R\Box B)\ :\ \ B\to B^{\prime}\nn
&&t^{-1}\colon =(L^r\stac{R} \cup_S)\ci \sqcap(L^r\stac{R} S,S^r)\ci
(\cap_L\ci \eta_R\Box B^{\prime}):\ B^{\prime}\to B.
\nonumber\eea
What is more $t:L^r\to S^r$ is a right ${\bf R}$-module isomorphism. This
means that we can choose the objects $B$ and $B^{\prime}$ to be the same and
require the equivalent conditions on $\M$-morphisms
\bea &&t=B\nn
\Leftrightarrow&&\cup_L\ci \sqcap(L,L^r)=\cup_S\ci \sqcap(S,S^r)\nn
\Leftrightarrow&&\cap_L\ci \eta_R=\cap_S\ci \eta_R
\nonumber\eea
to hold true. (Notice that then the $\M$-objects $S^r\stac{R} S$, $S^r\stac{R}
L$,  $L^r\stac{R} S$ and $L^r\stac{R} L$ are all equal.) The identity
$\lambda_{S^r}= \mu_B\ci \Sigma_{B,B}\ci (q\Box B)$ is checked by direct
calculation setting $q\colon =\lambda_{S^r}\ci(R\Box \eta_B)$ and using the
forms
\bea && \lambda_{S^r}=(S^r\stac{R} \cup_S)\ci \sqcap(S^r\stac{R} S,S^r)\ci
[\rho_{S^r\ot S}\ci(\cap_S\ci \eta_R\Box R)\Box B]\nn
&&\mu_B=(L^r\stac{R}\cup_L)\ci(L^r\stac{R} L\stac{R} \cup_L\stac{R} L^r)\ci
[(L^r\stac{R}\gamma)\ci \cap_L\stac{R}\sqcap(L^r,L^r)]\nn
&&\eta_B=(L^r\stac{R}\pi)\ci\cap_L\ci \eta_R
\nonumber\eea 
and the identity $\rho_{M_A}\ci (\gamma\Box s\Box \eta_A)=\gamma\ci \rho_S$.

Since $\pi\ci s=R$ also $s^r\ci \pi^r=R$ which proves that $s^r$ is an
epimorphism.
In order to prove the identity (\ref{psir}) let us introduce the ${\bf
B}$-${\bf R}$ bimodule 
$$ {\bf P}=\left(L^r\stac{R} S^r, \lambda_{N_A}\ci(\eta_B\Box B\Box
(L^r\stac{R} S^r)), \rho_{L^r\ot S^r}\right).$$
The identity $\psi\ci \sqcap(S,L)=\rho_{M_A}\ci(\gamma\Box A\Box \eta_A)$
implies that 
\bea&& (S^r\stac{R}\mu_S)\ci \sqcap(S^r\stac{R}S,S)\ci [(\gamma^r\stac{R} S)\ci
(L^r\stac{R}\cap_L)\Box A]\ci \Sigma_{A,B}=\nn
&&\qquad\qquad\qquad\qquad
(S^r\stac{R} L\stac{R}\cup_S)\ci \sqcap(S^r\stac{R}L\stac{R}S,S^r)\ci
[(S^r\stac{R}\psi)\ci \cap_S\stac{R} L)\Box B],
\nonumber\eea
from which it follows that
\bea&&(L^r\stac{R}S^r\stac{R}L\stac{R}\cup_S)\ci 
\sqcap(L^r\stac{R}S^r\stac{R}L\stac{R}S,S^r)\ci
[(L^r\stac{R}S^r\stac{R}\psi)\ci(L^r\stac{R}\cap_S\stac{R}L)\ci\cap_L\ci 
\eta_R\Box B]\nn
&&\qquad\qquad\qquad\qquad
=(L^r\stac{R}S^r\stac{R}\mu_S)\ci \lambda_{P\ot S\ot S}\ci [B\Box
  (S^r\stac{R}\cap_S\stac{R} S)\ci \cap_S\ci \eta_R],
\nonumber\eea
hence (\ref{psir}).

Applying the `right handed' version of Theorem \ref{entw-bgd} the properties
of the right entwining structure $({\bf L}^r,{\bf S}^r,\psi^r)$ proven so far
are equivalent to the statement that $({\bf B},{\bf R},\pi^r,q,\mu^r,s^r)$ is
a right bialgebroid in $\M$.
\hfill$\Box$

\bigskip

{\bf Acknowledgment.} We thank T. Brzezi\'nski and K. Szlach\'anyi for
their promoting questions.

The author was supported by the Hungarian Scientific
Research Fund OTKA -- T 020 285, FKFP -- 0043/2001  and the Bolyai
J\'anos Fellowship.

\end{document}